\documentclass[leqno]{amsart}
\usepackage{amssymb,color}
\usepackage{amssymb, amsmath}
\usepackage[final]{graphicx}
\usepackage{stackrel}
\makeatletter
\def\subsection{\@startsection{subsection}{2}%
\z@{.5\linespacing\@plus.7\linespacing}{.5\linespacing}%
  {\normalfont\bfseries}}
\makeatother
\usepackage[utf8]{inputenc}
\numberwithin{equation}{section}
\usepackage{color}

\newtheorem{theorem}{Theorem}

\newtheorem{remark}[theorem]{Remark}
\newtheorem{definition}[theorem]{Definition}
\newtheorem{corollary}[theorem]{Corollary}

\begin{document}
\title{Self-dual weakly-Einstein manifolds}
\author[ ]{JeongHyeong Park and Rodrigo Mari\~no Villar}
\address{JeongHyeong Park\\ Department of Mathematics, Sungkyunkwan University, Suwon 16419, Korea}
\email{parkj@skku.edu}
\address{Rodrigo Mari\~no Villar\\ Faculty of teaching training, University of Santiago de Compostela, Ramón Ferreiro Av., 27001, Lugo, Spain}
\email{rodrigo.marino.villar@usc.es}
\begin{abstract}

We describe the pointwise curvature of self-dual weakly Einstein four-manifolds through the existing algebraic classification of weakly Einstein curvature tensors. We identify generalized Singer-Thorpe (GST) bases with the simultaneous diagonalizing bases used in that classification and express the resulting Weyl spectra in GST parameters.The resulting GST statements cover all Ricci eigenvalue multiplicities and remain valid at Einstein points. We also recall Besse's classical harmonic Weyl dichotomy and its consequence in the weakly Einstein setting: a connected self-dual weakly Einstein four-manifold with harmonic Weyl curvature is either Einstein or locally conformally flat and scalar-flat.

\end{abstract}
\maketitle
\section{Introduction}

On an $n$-dimensional Riemannian manifold $(M,g)$, consider the
pointwise quadratic curvature condition

\begin{equation}\label{WE}
\check R=\frac{\|R\|^2}{n}g,
\end{equation}
where $\check R(X,Y)=R_{Xijk}R_Y^{ijk}$.
In dimension four, we have the curvature identity

\begin{equation}\label{eq1}
\left(\check R-\frac{\|R\|^2}{4}g\right)
+\tau\left(\rho-\frac{\tau}{4}g\right)
-2\left(\check\rho-\frac{\|\rho\|^2}{4}g\right)
-2\left(R[\rho]-\frac{\|\rho\|^2}{4}g\right)=0,
\end{equation}
where $\check\rho_{ij}=\rho_{ia}\rho^a_j$ and $(R[\rho])_{ij}=R_{iabj}\rho^{ab}$.
This identity comes from the Chern--Gauss--Bonnet formula \cite{EPS}.
Consequently, every Einstein metric in dimension four satisfies \eqref{WE}.

Motivated by curvature identities in dimension four, Euh, Park, and
Sekigawa introduced the following notion and exhibited non-Einstein
examples \cite{EPS013}.
\begin{definition}
A Riemannian manifold $(M,g)$ of dimension $n$ is called weakly-Einstein if
$\check R=\dfrac{\|R\|^2}{n}g$.
\end{definition}
This definition is pointwise and does not require compactness or
constancy of $\|R\|^2$.
An Einstein metric satisfying \eqref{WE} with constant $\|R\|^2$
is called super-Einstein
\cite{GW}.


For metrics with parallel Ricci tensor, condition~\eqref{WE} has the
following variational characterization. Let $g$ be a unit-volume metric
on a compact connected $n$-manifold without boundary, and suppose that
$\nabla\rho=0$. By \cite[Corollary~4.72(i)]{Besse}, $g$ is a critical
point of $\mathcal F(g)=\int_M |R|^2\,dV_g$ on the space of all
unit-volume metrics if and only if condition~\eqref{WE} holds.

In dimension four, \cite[(2.7)]{DKP_MJM} shows, without assuming parallel
Ricci curvature, that no proper weakly Einstein metric on a compact
connected manifold without boundary is a critical point of $\mathcal F$.
Here a weakly Einstein metric is called proper if it is neither Einstein
nor locally conformally flat with vanishing scalar curvature.
Since $\mathcal F$ is scale invariant in dimension four, this conclusion
also applies to critical points under the unit-volume constraint.

Garc\'ia-R\'io, Haji-Badali, Mari\~no-Villar and V\'azquez-Abal
proved that a non-Einstein locally conformally flat four-dimensional Riemannian
manifold is weakly Einstein if and only if its scalar curvature vanishes
\cite{EARE}.


A natural question is how the weakly Einstein condition constrains the
remaining Weyl component under half-conformal flatness, and when these
constraints force full conformal flatness.

A rigidity result is already available in the K\"ahler setting. Derdzinski,
Euh, Kim and Park proved in \cite[Theorem~1.2]{DEKP} that every self-dual
weakly Einstein K\"ahler surface, with respect to its complex orientation, is
locally symmetric. Up to local homothety, such a surface is the standard
complex projective or complex hyperbolic plane, or a product of two real
surfaces of opposite constant Gaussian curvatures, including the flat case.

In this paper, we study half-conformally flat weakly Einstein four-manifolds
without the K\"ahler assumption. Reversing the orientation if necessary,
we work with $W^-=0$.
The algebraic classification of weakly Einstein curvature tensors in
\cite[Theorems~1.3--1.5]{DerPark} already determines the pointwise
possibilities for the self-dual Weyl spectrum at non-Einstein points.
For a weakly Einstein curvature tensor, GST bases are precisely the
simultaneous diagonalizing bases of \cite[Theorem~1.1]{DerPark};
we verify this equivalence in Section~\ref{sec:gst}.
The bases used in the classification additionally arrange the eigenvalues
and parameters into the stated normal forms. The GST sign patterns and
the three classification cases organize the possibilities differently.
The GST statements in Section~\ref{sec:gst} are derived from this
classification and are recorded as corollaries.
{We follow the curvature-convention used in \cite{EPS, EPS013}, and  make the GST parameter restrictions explicit.}
The GST formulas also cover degenerate Ricci spectra and Einstein points.


We collect the self-dual consequences of the algebraic classification and
the classical harmonic Weyl dichotomy in the following theorem.
Parts~\textup{(i)} and~\textup{(ii)} are specializations of
\cite[Theorems~1.3--1.5]{DerPark}, {and} part~\textup{(iii)} combines
\cite[Corollary~16.29]{Besse} with \cite[Theorem~2(i)]{EARE}.
The spectral formulas are derived in Section~\ref{sec:gst}, and the
proof is completed in Section~\ref{sec:harmonic}.

\begin{theorem}\label{thm:main}\label{thm:sd-harmonic-weyl}
Let $(M^4,g)$ be a connected, oriented, self-dual weakly Einstein
Riemannian manifold. {Let $\tau$ be its scalar curvature and
$e=\rho-\tfrac14\tau g$ be its trace-free Ricci tensor.}
The spectra below are unordered.
\begin{enumerate}
\item[\textup{(i)}]
At any point $p$ with $e(p)\ne0$ at which the Ricci endomorphism does
not have two distinct eigenvalues of multiplicity two,
$$
\operatorname{spec}(W^+(p))=
\left(-\frac{\tau(p)}3,\frac{\tau(p)}6,\frac{\tau(p)}6\right).
$$
In particular, $W(p)=0$ if and only if $\tau(p)=0$.
\item[\textup{(ii)}]
At any point $p$ with $e(p)\ne0$ at which the Ricci endomorphism has
two distinct eigenvalues of multiplicity two, there exists
$\xi=\xi(p)\in\mathbb R$ such that
$$
\operatorname{spec}(W^+(p))=
\left(\frac{\tau(p)}6,\frac{\tau(p)}6-2\xi,
      2\xi-\frac{\tau(p)}3\right).
$$
Here, $W(p)=0$ if and only if $\tau(p)=\xi(p)=0$.
\item[\textup{(iii)}]
If $\delta W^+=0$, then either $g$ is Einstein, or $g$ is locally
conformally flat and scalar-flat.
\end{enumerate}
\end{theorem}

\section{Preliminaries}
Let $(M^4,g)$ be an oriented Riemannian manifold, and let $\nabla$ be its Levi-Civita connection and $R$ be its curvature tensor $R(X,Y)=[\nabla_X,\nabla_Y]-\nabla_{[X,Y]}$.
The Ricci tensor $\rho$, scalar curvature $\tau$, and Ricci operator $Q$ are defined by $$\rho(X,Y) = \operatorname{tr}\{Z\longmapsto R(Z,X)Y\}, \qquad \tau=\operatorname{tr}_g\rho, \qquad \rho(X,Y)=g(QX,Y). $$ For an orthonormal frame $\{e_1,e_2,e_3,e_4\}$, we write $ R_{ijkl}=g(R(e_i,e_j)e_k,e_l). $ At every point $p\in M$, the metric and orientation determine the Hodge operator $$ *:\Lambda^2T_p^*M\longrightarrow\Lambda^2T_p^*M. $$ Since $*^2=\operatorname{Id}$ on two-forms in dimension four, we have the orthogonal decomposition $$ \Lambda^2T_p^*M=\Lambda_p^+\oplus\Lambda_p^-, \qquad \Lambda_p^\pm = \{\eta\in\Lambda^2T_p^*M:* \eta=\pm\eta\}. $$ We regard the Weyl tensor as a self-adjoint endomorphism of $\Lambda^2T_p^*M$. {It} preserves the above decomposition, and hence $$ W=W^+\oplus W^-, \qquad W^\pm:\Lambda_p^\pm\longrightarrow\Lambda_p^\pm.$$ The metric $g$ is locally conformally flat if and only if $W=0$, and it is called half-conformally flat if $$W^+=0 \qquad\text{or}\qquad W^-=0. $$ With respect to the chosen orientation, the conditions $W^-=0$ and $W^+=0$ are called self-duality and anti-self-duality, respectively.

{
\section{Algebraic consequences of the classification}\label{sec:gst}
We choose the orientation so that $W^-=0$.
The pointwise results in this section are consequences of the algebraic
classification in \cite[Theorems~1.3--1.5]{DerPark}. We first specialize that classification
to the self-dual case and then express its consequences in GST notation.
All arguments are pointwise, and the GST statements below also cover points
where $e=\rho-\tfrac14\tau g$ vanishes.

\begin{definition}[\cite{EPS013}]
An orthonormal basis $\{e_1,e_2,e_3,e_4\}$ is a generalized
Singer--Thorpe (GST) basis if
\begin{equation}\label{eq:gst-definition}
R_{ijjk}=0\quad(i\ne k),\qquad
R_{1212}^2=R_{3434}^2,\quad
R_{1313}^2=R_{2424}^2,\quad
R_{1414}^2=R_{2323}^2.
\end{equation}
\end{definition}
The existence of a GST basis is equivalent to the weakly Einstein condition
\cite{EPS013}. This is a pointwise algebraic characterization, which we also
apply to auxiliary algebraic curvature tensors below. A GST basis is a
Ricci eigenbasis, and the induced bases of $\Lambda^\pm$ diagonalize $W^\pm$.
For a positive GST basis, {we} write
\begin{equation}\label{eq:gst-bivector-bases}
 E_1^\pm=e^1\wedge e^2\pm e^3\wedge e^4,\quad
 E_2^\pm=e^1\wedge e^3\mp e^2\wedge e^4,\quad
 E_3^\pm=e^1\wedge e^4\pm e^2\wedge e^3.
\end{equation}
These vectors have squared norm $2$. Consequently, the diagonal entries of
$W^\pm_{ij}=\langle W^\pm E_i^\pm,E_j^\pm\rangle$ are twice the
corresponding operator eigenvalues.

\subsection{Equivalence of adapted bases}\label{subsec:adapted-bases}

For any weakly Einstein algebraic curvature tensor, a positive
orthonormal basis is a GST basis if and only if it diagonalizes $e$ and
its associated bases $E_i^\pm$ diagonalize both $W^+$ and $W^-$.
Indeed, the GST condition $R_{ijjk}=0$ for $i\ne k$ implies that the
Ricci tensor is diagonal and that all Weyl components with exactly three
distinct indices vanish. The induced matrices of $W^\pm$ are therefore
diagonal. Conversely, simultaneous diagonalization gives vanishing of
the three-distinct-index Weyl components. Since $e$ is diagonal, the
curvature decomposition then gives $R_{ijjk}=0$ for $i\ne k$.
For a tensor with this vanishing property, direct contraction gives
\begin{equation}\label{eq:gst-contraction}
4\bigl(R_{1212}^2-R_{3434}^2\bigr)
=\check R_{11}+\check R_{22}-\check R_{33}-\check R_{44}=0,
\end{equation}
where the last equality uses the weakly Einstein condition.
The analogous identities obtained by permuting the indices give
$R_{1313}^2=R_{2424}^2$ and $R_{1414}^2=R_{2323}^2$.
Thus the bases of \cite[Theorem~1.1]{DerPark} are precisely GST bases,
and those in \cite[Theorems~1.3--1.5]{DerPark} additionally put the
eigenvalues and parameters into normal form. This equivalence does not
require half-conformal flatness and remains valid at Einstein points.
The different case divisions therefore concern the GST sign patterns
and the classification parameters, rather than different kinds of bases.

\subsection{Self-dual specialization of the classification}\label{subsec:self-dual-classification}
The conventions of \cite{DerPark}, denoted by hats, are related to ours by
\begin{equation}\label{eq:curvature-convention}
\widehat R=-R,\qquad \widehat\rho=\rho,\qquad
s=\tau,\qquad \widehat W=-W.
\end{equation}
At a point with $e\ne0$, distinguish case (a), in which $e$ does not have
two distinct double eigenvalues, and case (b), in which it does.
Then \cite[Theorems~1.3--1.5]{DerPark} gives
\begin{equation}\label{eq:classified-spectra}
\operatorname{spec}(W^+)=
\begin{cases}
 (\tau/6,\tau/6,-\tau/3),
 &\text{in case (a)},\\
 (\tau/6,\tau/6-2\xi,2\xi-\tau/3),
 &\text{in case (b)},
\end{cases}
\end{equation}
where $\xi\in\mathbb R$ and the spectra are unordered.
Indeed, the eigenvalues of $\widehat W^\pm$ in the normal forms have the
form $a_j\pm c_j$, $j=2,3,4$. In the first case,
$$
(a_2,a_3,a_4)=(-\tau/12,-\tau/12,\tau/6),
$$
by Theorem~1.3 when $\tau=0$ and Theorem~1.4 when $\tau\ne0$.
In the second case, Theorem~1.5 gives
$$
(a_2,a_3,a_4)=(-\tau/12,\xi-\tau/12,-\xi+\tau/6).
$$
The condition $\widehat W^-=0$ forces $c_j=a_j$, so the eigenvalues of
$W^+$ are $-2a_j$, proving \eqref{eq:classified-spectra}.
In particular, in the first case $W=0$ if and only if $\tau=0$;
in the second case $W=0$ if and only if $\tau=\xi=0$.

\subsection{GST parameters and degenerate Ricci spectra}\label{subsec:gst-parameters}
{We set}
$$
 R_{1212}=\alpha,\qquad R_{1313}=\beta,\qquad R_{1414}=\gamma,
$$
and {we write} $R_{3434}=\epsilon_{12}\alpha$,
$R_{2424}=\epsilon_{13}\beta$, and
$R_{2323}=\epsilon_{14}\gamma$, with $\epsilon_{1j}\in\{1,-1\}$.
The pattern $(1,1,1)$ has $e=0$. The other patterns, up to permutation,
are $(-1,-1,-1)$, $(1,-1,-1)$ and $(1,1,-1)$; these patterns need not
be disjoint when some curvature components vanish.
The Ricci contraction gives
\begin{equation}\label{eq:gst-ricci}
\begin{aligned}
\rho_{11}&=-\alpha-\beta-\gamma,&
\rho_{22}&=-\alpha-\epsilon_{13}\beta-\epsilon_{14}\gamma,\\
\rho_{33}&=-\epsilon_{12}\alpha-\beta-\epsilon_{14}\gamma,&
\rho_{44}&=-\epsilon_{12}\alpha-\epsilon_{13}\beta-\gamma,\\
\tau&=-2\{\alpha(1+\epsilon_{12})+\beta(1+\epsilon_{13})
    +\gamma(1+\epsilon_{14})\}.
\end{aligned}
\end{equation}
For complementary pairs, the Weyl decomposition yields
\begin{equation}\label{eq:gst-complement}
W_{ijij}=\frac12(R_{ijij}+R_{klkl})+\frac{\tau}{12},
\qquad \{i,j,k,l\}=\{1,2,3,4\}.
\end{equation}
This agrees with the equality of complementary Weyl components in
\cite[(2.4)]{DerPark}.
Since $W^-=0$, the eigenvalues of $W^+$ in the basis of $\Lambda^+$ induced by the GST basis are $2W_{1212}$, $2W_{1313}$, and $2W_{1414}$.
Substituting the corresponding GST sign patterns into \eqref{eq:gst-ricci} and \eqref{eq:gst-complement} also gives a direct derivation of the vanishing of $W$ in Corollary 4 and the spectra in Corollaries 5 and 7.
If $\nu_i=e(e_i,e_i)$, the same decomposition reads
\begin{equation}\label{eq:replace-ricci}
R_{ijij}=W_{ijij}-\frac12(\nu_i+\nu_j)-\frac{\tau}{12}.
\end{equation}

The converse parts of the three classification theorems are organized by
Ricci multiplicities, rather than by GST sign patterns. In particular,
Theorems~1.3(c) and~1.4(c) of \cite{DerPark} do not cover two distinct
double eigenvalues, and the non-Einstein classification does not apply
at $e=0$. The following elementary algebraic reduction allows us to
apply them without imposing any extra geometric hypothesis.

Fix a GST basis and keep $W$ and $\tau$ unchanged in the curvature
decomposition, while replacing $e$ by another diagonal trace-free tensor.
Equation~\eqref{eq:replace-ricci} shows that all complementary sectional
sums remain unchanged, as do the vanishing components with three distinct
indices. For type $(-1,-1,-1)$, any such replacement preserves the GST
pattern. For type $(1,-1,-1)$, the replacement
$$
\widetilde e=\operatorname{diag}(-3,3,-1,1)
$$
also preserves the equal pair $R_{1212}=R_{3434}$, since
$\widetilde\nu_1+\widetilde\nu_2=
 \widetilde\nu_3+\widetilde\nu_4=0$.
Thus this replacement works for both types and produces a weakly Einstein
algebraic curvature tensor with four distinct Ricci eigenvalues, with
exactly the original $W$ and $\tau$.
For type $(1,1,-1)$, the replacement
$\widetilde e=\operatorname{diag}(-r,r,r,-r)$, $r>0$, preserves both
equal complementary pairs and the opposite pair. {This} produces two distinct
double eigenvalues, again with $W$ and $\tau$ unchanged.
These are replacements of algebraic curvature tensors at a single point;
no perturbation of the Riemannian metric, geometric realization, or density
of a generic set is asserted or needed.

\begin{corollary}\label{thm:three-minus}
Let $(M^4,g)$ be a non-Einstein oriented Riemannian manifold admitting
at every point a GST basis of type $(-1,-1,-1)$.
Then $(M,g)$ is half-conformally flat and weakly Einstein if and only if
it is locally conformally flat and scalar-flat.
\end{corollary}
\begin{proof}
Suppose first that $g$ is half-conformally flat and choose the self-dual
orientation. The GST relations give $\tau=0$ by \eqref{eq:gst-ricci}.
At any point, make the algebraic replacement described above, preserving
$W$ and $\tau$. The auxiliary tensor has four distinct Ricci eigenvalues,
so \cite[Theorem~1.3(c)]{DerPark} applies. In that normal form the spectra
of $\widehat W^+$ and $\widehat W^-$ are mutually opposite.
Since $\widehat W^-=0$, both vanish. Hence the original $W$ also vanishes.
Conversely, $W=0$ and $\tau=0$ imply the weakly Einstein condition by
\cite[(1.1)]{DerPark}, and half-conformal flatness is immediate.
\end{proof}

\begin{corollary}\label{thm:two-minus}
Let $(M^4,g)$ be a non-Einstein oriented, half-conformally flat Riemannian
manifold admitting at every point a GST basis of type $(1,-1,-1)$. Then, with the self-dual orientation,
$$
\operatorname{spec}(W^+)=(-\tau/3,\tau/6,\tau/6).
$$
In particular, $(M,g)$ is locally conformally flat if and only if it is
scalar-flat.
\end{corollary}
\begin{proof}
Here $\tau=-4\alpha$. Preserve $W$ and $\tau$ and replace $e$ as above
so that the Ricci eigenvalues are distinct.
If $\tau=0$, \cite[Theorem~1.3(c)]{DerPark} and self-duality give $W=0$.
If $\tau\ne0$, \cite[Theorem~1.4(c)]{DerPark}, specialized as in
\eqref{eq:classified-spectra}, gives the displayed spectrum.
Since the Weyl tensor was unchanged, the formula holds for the original
tensor at every point. The last assertion follows from $W^-=0$.
\end{proof}

\begin{remark}[Relation to the exceptional algebraic family]\label{rem:normal}
The preceding two corollaries are consequences of the classification,
including at Einstein points and at points with two distinct double
Ricci eigenvalues. These cases are covered by the algebraic reduction
above.
For type $(1,-1,-1)$, the trace-free Ricci eigenvalues are
\begin{equation}\label{eq:e-two}
(-\beta-\gamma,\beta+\gamma,-\beta+\gamma,\beta-\gamma).
\end{equation}
Thus $e=0$ precisely when $\beta=\gamma=0$. At a non-Einstein point,
setting $\lambda=|\beta|+|\gamma|$ and
$\mu=\bigl||\beta|-|\gamma|\bigr|$ puts these eigenvalues in the
nondecreasing order $(-\lambda,-\mu,\mu,\lambda)$.
If $\beta\gamma\ne0$, then $\lambda>\mu$ and
\cite[Lemma~2.1(a)]{DerPark} gives the spectrum directly.
If $\beta\gamma=0$, then $\lambda=\mu>0$ and part~(b) initially permits
an extra parameter. The GST condition restricts this parameter to the
spectrum of Corollary~\ref{thm:two-minus}, as the algebraic reduction
shows. Equivalently, in the parameterization of
\eqref{eq:classified-spectra}, one may take $\xi=0$ or $\xi=\tau/4$,
depending on the ordering of the two mixed planes.
In particular, the exceptional family is not omitted from the argument.
\end{remark}

\begin{corollary}\label{thm:one-minus}
Let $(M^4,g)$ be a non-Einstein oriented, half-conformally flat Riemannian
manifold admitting at every point a GST basis of type $(1,1,-1)$. In such a basis, {let}
$R_{1212}=\alpha$, $R_{1313}=\beta$, and $R_{1414}=\gamma$.
Then
\begin{equation}\label{eq:one-minus-spectrum}
\tau=-4(\alpha+\beta),\qquad
\operatorname{spec}(W^+)=
(\tau/6+2\alpha,\tau/6+2\beta,\tau/6).
\end{equation}
Consequently, $(M,g)$ is locally conformally flat if and only if it is
scalar-flat and $\alpha=\beta$, or equivalently if and only if
$\alpha=\beta=0$ at every point.
\end{corollary}
\begin{proof}
Equation~\eqref{eq:gst-ricci} gives
\begin{equation}\label{eq:e-one}
\begin{aligned}
\rho_{11}=\rho_{44}&=-\alpha-\beta-\gamma,\\
\rho_{22}=\rho_{33}&=-\alpha-\beta+\gamma,\\
\tau&=-4(\alpha+\beta),
\end{aligned}
\qquad e=\operatorname{diag}(-\gamma,\gamma,\gamma,-\gamma).
\end{equation}
For $\gamma\ne0$, an orientation-preserving reordering of the GST basis
puts the negative eigenvalues of $e$ first. With the mixed planes ordered
appropriately, \eqref{eq:gst-complement} gives, in the convention of
\cite{DerPark},
$$
(a_2,a_3,a_4)=(-\tau/12,-\alpha-\tau/12,-\beta-\tau/12).
$$
Thus Theorem~1.5 of \cite{DerPark} has parameter $\xi=-\alpha$,
with $\beta=\xi-\tau/4$.
Its self-dual specialization \eqref{eq:classified-spectra} yields
\eqref{eq:one-minus-spectrum}, up to permutation.
When $\gamma=0$, replace $e$ by
$\operatorname{diag}(-r,r,r,-r)$, $r>0$, as above.
This leaves $W$, $\tau$, $\alpha$ and $\beta$ unchanged, so the same
argument applies. Finally, the spectrum vanishes exactly when
$\tau=\alpha=\beta=0$. Since $\tau=-4(\alpha+\beta)$, this is equivalent
to the two conditions stated in the corollary.
\end{proof}

The corresponding unnormalized component matrix is
\begin{equation}\label{eq:one-minus-matrix}
(W^+_{ij})=\frac43\operatorname{diag}
(2\alpha-\beta,2\beta-\alpha,-\alpha-\beta)
=\frac13\operatorname{diag}(\tau+12\alpha,\tau+12\beta,\tau).
\end{equation}
At $\tau=0$, the exceptional algebraic spectrum is
$(0,2\xi,-2\xi)$, up to ordering. Hence scalar-flatness alone does not
force the Weyl tensor of an algebraic self-dual weakly Einstein tensor
to vanish. The classification supplies algebraic possibilities; it does
not by itself establish their realization by metrics on an open set.
We next recall the classical restriction under harmonic Weyl curvature.
}

{
\section{Self-duality and harmonic Weyl curvature}\label{sec:harmonic}
We choose the orientation so that $W^-=0$, and {we write}
$e=\rho-\tfrac14\tau g$ for the trace-free Ricci tensor.
The dichotomy in \cite[Corollary~16.29]{Besse} applies without the
weakly Einstein assumption; the latter is needed only for the
scalar-flatness conclusion below.

\begin{proof}[Proof of Theorem~\ref{thm:main}]
The spectral assertions in parts~\textup{(i)} and~\textup{(ii)} follow
from \eqref{eq:classified-spectra}. Since $W^-=0$, the tensor $W(p)$
vanishes exactly when all three eigenvalues of $W^+(p)$ vanish.
The displayed triples therefore give the stated vanishing criteria.

For part~\textup{(iii)}, the hypothesis $\delta W^+=0$ gives
$\delta W=0$ because $W^-=0$. By \cite[Corollary~16.29]{Besse},
$g$ is Einstein or locally conformally flat. If $g$ is non-Einstein,
the weakly Einstein assumption and \cite[Theorem~2(i)]{EARE} give
$\tau\equiv0$.
\end{proof}

The Codazzi-tensor argument underlying the dichotomy is classical as well.
The condition $\delta W=0$ makes the Schouten tensor
$$
A=\frac12\left(\rho-\frac{\tau}{6}g\right)
$$
Codazzi. Its trace-free part is $e/2$, so it is not proportional to $g$
on $U=\{e\ne0\}$. Applying \cite[Theorem~16.21(ii)]{Besse} on $U$
gives equality of the spectra of $W^+$ and $W^-$ there; see also
\cite[Corollary~1]{Derdzinski1988}. Self-duality therefore gives $W=0$
on $U$. If $g$ is non-Einstein, $U$ is a nonempty open set.
Unique continuation for the harmonic Weyl tensor extends this vanishing
to the connected manifold $M$. This global step is already included in
\cite[Corollary~16.29 and the proof of Theorem~16.31]{Besse}.
{
\begin{remark}[Open-set realization of GST types]\label{rem:sd-gst}
Let $O$ be an open set on which $W^-=0$, $e\ne0$, and $W\ne0$.
Then $\{p\in O:(\delta W^+)(p)\ne0\}$ is open and dense in $O$.
Indeed, if $\delta W^+$ vanished on a nonempty open subset of $O$,
the Codazzi-tensor argument recalled above would give $W=0$ there.
This applies to both GST types $(1,-1,-1)$ and $(1,1,-1)$.

For type $(1,-1,-1)$, Corollary~\ref{thm:two-minus} gives the
eigenvalues of $W^+$ as $(-\tau/3,\tau/6,\tau/6)$. Thus $\tau\ne0$
on $O$. Let $P$ be the orthogonal projection onto the simple
eigenline and $I$ the identity on $\Lambda^+$. Then
$$
   W^+=\frac{\tau}{6}(I-3P),\qquad
   \nabla W^+=\frac16\,d\tau\otimes(I-3P)-\frac{\tau}{2}\nabla P.
$$
If $\tau$ is constant on $O$, then $\nabla P\ne0$ wherever
$\delta W^+\ne0$, hence on an open dense subset of $O$.
Whether a self-dual weakly Einstein metric can realize this GST type
throughout such an open set without the harmonicity assumption remains
a separate question.
\end{remark}}

}

\section*{Acknowledgements}
The authors thank Andrzej Derdziński for pointing out how Besse’s Theorem 16.21(ii) applies to Theorem 2(iii) and how certain curvature formulas also follow from the algebraic classification. J. H. Park was supported by the National Research Foundation of Korea (NRF)
grant funded by the Korean government (MSIT)
(No.~RS-2024-00334956).


\begin{thebibliography}{99}
\bibitem{Besse}
A.~L.~Besse, \emph{Einstein Manifolds},
Ergebnisse der Mathematik und ihrer Grenzgebiete (3), vol.~10,
Springer-Verlag, Berlin, 1987.
\bibitem{DKP_MJM}
A.~Derdzinski, S.~Kim, and J.~H.~Park,
\emph{Scalar-flat K\"ahler surfaces whose Weyl tensor annihilates the Ricci form},
Mediterr. J. Math., to appear.
\bibitem{DerPark}
	{ A.~Derdzinski, J. H. ~Park and W.~Shin,
\emph{Weakly Einstein curvature tensors},
Glasg.  Math.  J., to appear.}

\bibitem{DEKP}
A.~Derdzinski, Y.~Euh, S.~Kim, and J.~H.~Park,
\emph{On weakly Einstein K\"ahler surfaces},
Int. J. Math. \textbf{36} (2025), no.~14, Article No.~2550057.
\bibitem{Derdzinski1988}
A.~Derdzi\'nski, \emph{Riemannian metrics with harmonic curvature on
$2$-sphere bundles over compact surfaces}, Bull. Soc. Math. France
\textbf{116} (1988), no.~2, 133--156.



\bibitem{EPS}
Y. Euh, J. Park, and K. Sekigawa,
\emph{A curvature identity on  a 4-dimensional Riemannian manifold},
{Results Math.} \textbf{63} (2013), 107--114.


\bibitem{EPS013}
Y. Euh, J. Park, and K. Sekigawa,
\emph{A generalization of a 4-dimensional Einstein manifold},
{Math. Slovaca} \textbf{63} (2013), 595--610.

\bibitem{EARE} E. Garc\'ia-R\'io, A. Haji-Badali, R. Mari\~no-Villar, M. E. V\'azquez-Abal, \emph{Locally conformally flat weakly-Einstein manifolds}.  Arch. Math. (Basel) \textbf{111} (2018), 549-559.

\bibitem{GW} A. Gray and T. J. Willmore, Mean-value theorems for Riemannian manifolds, Proc.
Roy. Soc. Edinburgh Sect. A \textbf{92} (1982), 343–364.
\end{thebibliography}
\end{document}